\documentclass{article}
\usepackage{spconf,amsmath,graphicx}
\RequirePackage{amsthm,amsmath,amsfonts}
\usepackage[T1]{fontenc}
\usepackage{color}
\usepackage{epsfig}
\usepackage[latin1]{inputenc}
\usepackage{rotating}
\usepackage{amssymb}
\usepackage{epsf}
\usepackage{graphicx}
\usepackage{psfrag}
\usepackage{subfigure}
\usepackage{url}
\usepackage{amsbsy}
\usepackage{mathtools}
\mathtoolsset{showonlyrefs,showmanualtags}

\newcommand{\ie}{\emph{i.e.}}

\newcommand{\cf}{\emph{c.f.}}
\newcommand{\sgn}{\mathrm{sgn}}

\renewcommand{\L}{{\mathcal L}}

\numberwithin{equation}{section}
\theoremstyle{plain}

\title{How to monitor and mitigate stair-casing in l1 trend filtering}

\name{Cristian R. Rojas and Bo Wahlberg\thanks{This work was partially supported by the Swedish Research Council and the Linnaeus Center ACCESS at KTH. The research leading to these results has received funding from The European Research Council under the European Community's  Seventh Framework program (FP7 2007-2013) / ERC Grant Agreement N. 267381.}}

\address{Department of Automatic Control and ACCESS Linnaeus Centre\\
School of Electrical Engineering, KTH Royal Institute of Technology, Sweden}
\begin{document}

\maketitle

\begin{abstract}
In this paper we study the estimation of changing trends in time-series using $\ell_1$ trend filtering. This method generalizes 1D Total Variation (TV) denoising for detection of step changes in means to detecting changes in trends, and it relies on a convex optimization problem for which there are very efficient numerical algorithms. It is known that TV denoising suffers from the so-called stair-case effect, which leads to detecting false change points. The objective of this paper is to show that $\ell_1$ trend filtering also suffers from a certain stair-case problem.  The analysis is based on an interpretation of the dual variables of the optimization problem in the method as integrated random walk. We discuss consistency conditions for $\ell_1$ trend filtering, how to monitor their fulfillment, and how to modify the algorithm to avoid the stair-case false detection problem.
\end{abstract}
\begin{keywords}
$\ell_1$ trend filtering, generalized lasso, TV denoising, Fused Lasso, change point detection.
\end{keywords}

\section{Introduction} \label{sec:intro}
We study the $\ell_1$ trend filtering method given (for $\lambda>0$) by
\begin{align} \label{eq:tv}
\min_{\{m_t\}_{t=1}^N} \frac{1}{2}\sum_{t=1}^N (y_t-m_t)^2+\lambda \sum_{t=3}^N |m_t - 2m_{t-1} + m_{t-2}|, \nonumber
\\ \end{align}
to estimate mean-trends in a time series dataset $\{y_t\}_{t=1}^N$ generated by the non-stationary Gaussian process
\begin{align}
y_t \sim \mathcal{N}(m_t, \sigma^2),
\end{align}
where the variance $\sigma^2 > 0$ is constant. It is assumed that the mean $\{m_t\}$ forms a piecewise linear sequence, \ie, a piecewise linear \emph{trend}. One way to measure the variability of a sequence $\{x_t, t=1, \dots, N\}$ is via its Total Variation (TV)\footnote{Another approach, for instance, is to specify the probability of a change and then use for example multiple model estimation methods \cite{gustafsson-2000a}.}:
\begin{align*}
\sum_{t=2}^N |x_t-x_{t-1}|.
\end{align*}
This is the $\ell_1$-norm of the first-difference sequence and can be seen as a convex approximation/relaxation of counting the number of changes. It is also directly related to measuring sparseness using the $l_1$-norm, as in the lasso method. Since the trend is assumed to be piecewise linear, we consider the second difference $x_t = m_t -2 m_{t-1}+m_{t-2}$, and we impose an $\ell_1$ penalty on $x_t$: $\sum_{t=3}^N |m_t - 2m_{t-1} + m_{t-2}|$. The fit to the data is measured by the least squares cost function
$\frac{1}{2}\sum_{t=1}^N (y_t-m_t)^2$,
which is related to the Maximum Likelihood (ML) cost function for the normally distributed case. The so-called $\ell_1$  trend filter \cite{Kim-Koh-Boyd-Gorinevsky-09} is given by minimizing a convex combination of these two cost functions, leading to (\ref{eq:tv}). This is a convex optimization problem with only one design parameter, namely $\lambda>0$. The TV cost will promote solutions for which $m_t - 2m_{t-1} + m_{t-2}=0$, \ie, a piecewise linear estimate (without jumps). As remarked in \cite{Kim-Koh-Boyd-Gorinevsky-09}, this method is related to the Hodrick-Prescott filter \cite{citeulike:11033787}, where an $\ell_2$ penalty on the second difference sequence is imposed; however, the $\ell_1$ norm is better at promoting sparsity, which translates here into a piecewise linear sequence $m_t$. The choice of the regularization parameter $\lambda$ is very important and provides a balance between fitting the data and stressing the structure constraint.  The same idea can be used for the multivariate case, \ie, for a vector valued stochastic process. The $\ell_1$ norm can then be replaced by a sum of norms, and
%
%
is known as \emph{sum-of-norms regularization} \cite{OhlssonLB:10}. For simplicity of presentation, however, we will focus on the univariate case. The $\ell_1$ trend filtering method is a special case of the generalized lasso method studied in \cite{Ryan11}. It is also related to spline approximations \cite{SDN06,Ryan14}.

The corresponding problem of detecting and estimating step-changes in means that are piecewise constant using
\begin{align}
\min_{m_1,\ldots, m_N} \frac{1}{2}\sum_{t=1}^N (y_t-m_t)^2+\lambda \sum_{t=2}^N |m_t - m_{t-1}|. \label{eq:mean-tv}
\end{align}
is more well studied. This method is called one-dimensional Total Variation (TV) denoising, Fused Lasso Signal Approximator or $l_1$ mean filtering \cite{citeulike:9394683,citeulike:9394684}. Some asymptotic convergence properties of the fused lasso are given in  \cite{Rinaldo2009}. In \cite{Rojas-Wahlberg-14b} it was rigorously shown that $l_1$ mean filtering  detection method  fails under well defined and intuitive conditions, namely  when two consecutive changes in the mean have the same sign (called a stair-case).  The objective of the current paper is to show that a similar problem also occurs for the $l_1$ trend filtering, but, even more importantly, how this problem can be monitored and mitigated. We propose an alternative  method to avoid this problem. The idea is to notice that the first and last detected change points in a sequence do not suffer from the stair-case effect. We therefore propose to restart the algorithm in a second step using only data in between these two detected change points, and then iteratively go through the whole sequence in the same way. This idea was inspired by \cite{Piotr14}, which uses random segmentation intervals.

In Section~\ref{sec:optim}, the optimality conditions for the method are derived, and Section~\ref{sec:interp} presents an interpretation of these conditions, based on which a consistency analysis is performed; for reasons of space, we only present a heuristic derivation, based on the analysis of a related problem (\cf, \cite{Rojas-Wahlberg-14b}), postponing the analytic details for a later publication. In Section~\ref{sec:monitoring} a modified scheme to remove fake change points is presented.
Section~\ref{sec:example} illustrates some examples of the method and its consistency, and Section~\ref{sec:conclusion} concludes the paper.

\section{Optimality Conditions} \label{sec:optim}

To derive the optimality conditions for the $\ell_1$ trend filter, we re-write \eqref{eq:tv} as
%
\begin{align} \label{eq:trend}
&\min\limits_{\{m_t\}_{t=1}^N,\{w_t\}_{t=2}^N} \; \displaystyle \frac{1}{2}\sum_{t=1}^N [y_t-m_t]^2+\lambda \sum_{t=3}^N |w_t|\\
&\qquad\quad \text{s.t.} \qquad\;\;  w_t=m_t - 2m_{t-1} + m_{t-2}, t = 3, \dots, N.
\end{align}

\vspace{-5mm}
\noindent To derive the optimality conditions, consider the Lagrangian

\vspace{-8mm}
\begin{multline*}
{\L}(\{m_t\}_{t=1}^N,\{w_t\}_{t=2}^N,\{z_t\}_{t=2}^{N-1}) 
=\frac{1}{2}\sum_{t=1}^N [y_t-m_t]^2+ \\
\lambda \sum_{t=3}^N |w_t| + \sum_{t=3}^N z_{t-1} (m_t - 2 m_{t-1}+m_{t-2} - w_t).
\end{multline*}
Minimizing $\L$ with respect to $m_1, \dots, m_N$, we obtain
\begin{align*}
&-(y_1 - m_1) + z_2 = 0, \\
&-(y_2 - m_2) - 2 z_2 + z_3 = 0, \\
&-(y_t - m_t) + z_{t-1} - 2 z_t + z_{t+1} = 0, \; t=3,\ldots, N-2,\\
&-(y_{N-1} - m_{N-1}) + z_{N-2} - 2 z_{N-1} = 0,  \\
&-(y_N - m_N) + z_{N-1} = 0.
\end{align*}
Iterating these equations backwards in $t$ gives
%
%
\begin{align} \label{eq:zt}
z_t &= \sum_{i=1}^{t-1} (t - i) [m_i-y_i], \quad t = 0, \dots, N+1,
\end{align}
with initial and end conditions $z_0 := z_1 := z_N := z_{N+1} := 0$. Thus, $\{z_t\}$ are a doubly integrated version of $\{m_t - y_t\}$, and correspond to the dual variables of the $\ell_1$ trend filtering method.

To minimize $\L$ with respect to $w_3, \dots, w_N$, we force the subgradient of $\L$ with respect to $w_t$ to equal $0$, which gives
\begin{align*}
z_{t-1} \left\{\begin{array}{ll}
= -\lambda, & w_t < 0,\\
\in [-\lambda, \lambda],& w_t = 0,\\
= \lambda, & w_t > 0.
\end{array}\right. \qquad t = 3, \dots, N,
\end{align*}
%
%
Therefore, since $w_t=m_t - 2m_{t-1} + m_{t-2}$,
\begin{align} \label{eq:optcondtrend}
&|z_t| \leqslant \lambda,\;\; t=2,\ldots , N-1, \nonumber\\
&|z_t| < \lambda\;\; \Rightarrow\: m_{t+1} - 2m_t + m_{t-1}=0, \\
&|z_t| = \lambda \;\; \Rightarrow\; \sgn(m_{t+1} - 2m_t + m_{t-1}) = \sgn(z_t). \nonumber
\end{align}
where $\sgn(x) := 1$ if $x >0$, $\sgn(x) := -1$ if $x < 0$ and $\sgn(0) := 0$. In the next section we will study the optimality conditions \eqref{eq:zt}, \eqref{eq:optcondtrend} in more detail, to derive consistency conditions.

\section{Interpretation and Consistency} \label{sec:interp}

\begin{figure}[t]
\centering
\includegraphics[width=0.7\columnwidth]{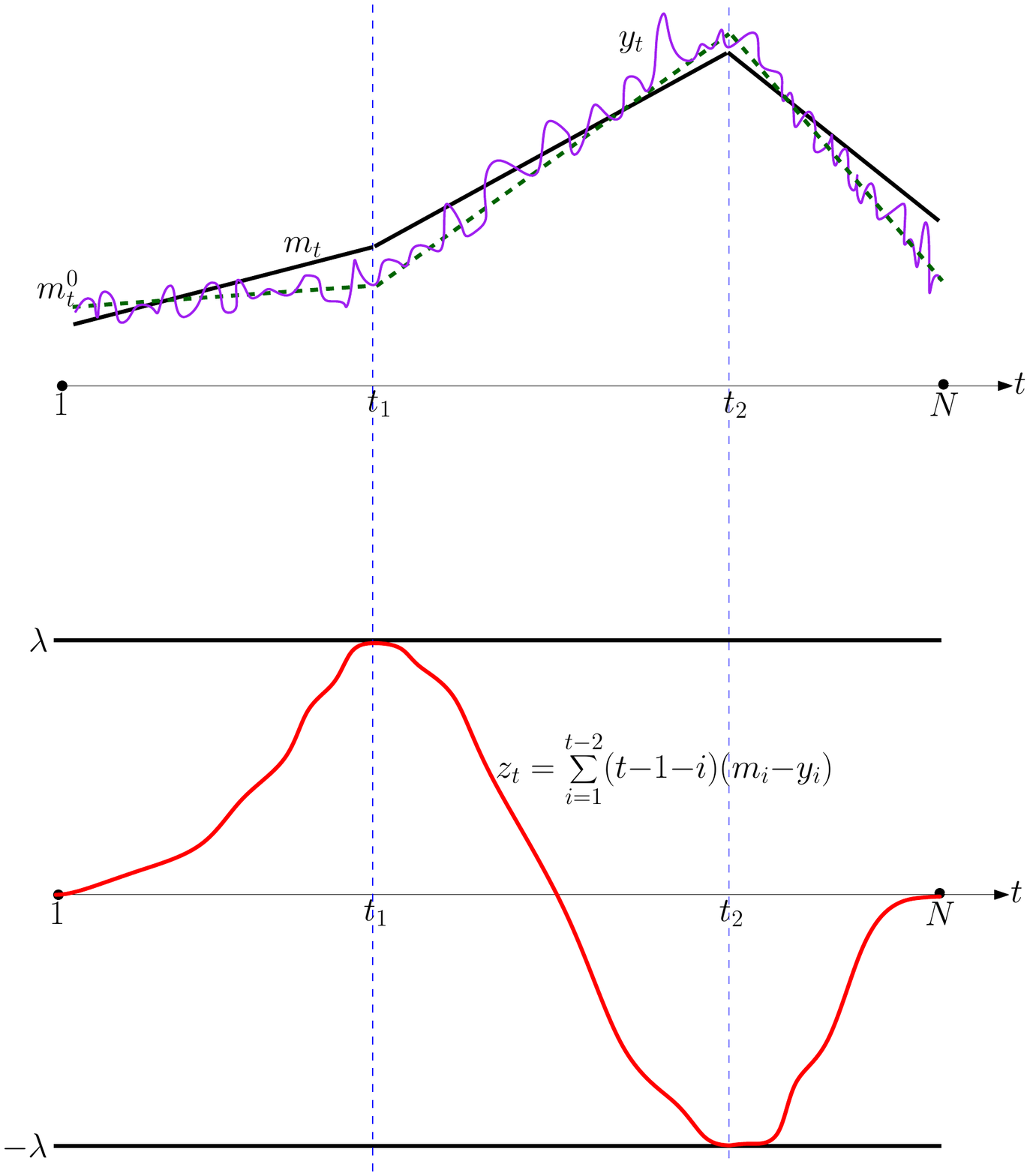}
\caption{Interpretation of optimality conditions (\ref{eq:zt}), (\ref{eq:optcondtrend}).}
\label{fig:1}
\end{figure}

The optimality conditions \eqref{eq:zt}, \eqref{eq:optcondtrend} can be interpreted according to the sketch of Fig.~\ref{fig:1}. From \eqref{eq:zt}, $z_t$ is essentially a doubly integrated version of $m_t - y_t$. If $m_t = m_t^o$, the true mean of $y_t$, then $z_t$ would be an integrated random walk process, since the term $m_t - y_t$ is essentially white Gaussian noise plus a deterministic term. In the general case, as $m_t$ and $m_t^o$ are both piecewise linear without jumps, the deterministic term is a discrete version of a cubic spline, \ie, a piecewise cubic polynomial with continuous derivatives of second order. Due to \eqref{eq:optcondtrend}, $z_t$ must always lie between $-\lambda$ and $\lambda$, and only touch the boundaries of this tube whenever there is a change in the slope of $m_t$; $z_t$ must equal $\lambda$ at $t_0$ if $m_{t_0+1} - m_{t_0} > m_{t_0} - m_{t_0-1}$ (\ie, if the slope of $m_t$ increases at $t_0 $), or $-\lambda$ if the reverse inequality holds. In addition, $z_0 = z_1 = z_N = z_{N+1} = 0$, which impose a series of interpolation constraints on the dual variables $z_t$. To satisfy these constraints, and those imposed by \eqref{eq:optcondtrend}, the estimate $m_t$ must suffer a bias whose integrated effect must be positive in segments where $z_t$ should go from $-\lambda$ (or $0$) to $\lambda$ (or $0$), and negative otherwise.

For $\lambda=0$, the method delivers $m_t = y_t$. As $\lambda$ is increased, the bias terms need to be increased so that $z_t$ at the change points can touch the boundaries $\pm \lambda$. However, as shown in Fig.~\ref{fig:1}, this leads to an estimated trend whose neighboring slopes at the change points differ less than the true slopes, and these differences decrease further as $\lambda$ is incremented, to the point where neighboring slopes coincide, and the neighboring segments are \emph{fused}. When $\lambda$ overcomes a prescribed value, called $\lambda_{\max}$, all segments are fused together, and $\ell_1$ trend filtering delivers a single linear trend for the entire dataset. This behavior resembles that of the fused lasso technique, as detailed in \cite{Rojas-Wahlberg-14b}.

To study the consistent recovery of the change points of $m_t^o$ (\ie, those values of $t$ for which $m_t^o - m_{t-1}^o$ changes\footnote{In many applications, it is important to know when the changes in trend have occurred. In addition, once the change points have been located, the trend can be consistently estimated by fitting a linear function to each individual data segment between the estimated change points.}), we consider the following asymptotic regime\footnote{These assumptions are made for simplicity, but they can be relaxed.}:
\begin{itemize}
\vspace{-2mm}
\item the number of samples $N$ tends to infinity;

\vspace{-2mm}
\item the variance $\sigma^2$ is kept constant (with respect to $N$);

\vspace{-2mm}
\item the number of change points $M$ is bounded; and

\vspace{-2mm}
\item the magnitude of the changes in slope of $m_t^o$ is bounded from above and from below.
\end{itemize}

Following the derivation in \cite{Rojas-Wahlberg-14b}, one can show that under these assumptions it is possible to recover the approximate location of all the change points \emph{if the consecutive changes in slope have all alternating signs}, as in Fig.~\ref{fig:1}. By ``approximate location'' we mean that, for a specific choice of $\lambda$, the estimated $m_t$ would have change points (perhaps more than $1$) at a distance $O(\epsilon)$ of the true change points of $m_t^o$, where $\epsilon$ is fixed but arbitrarily small, and no other estimated change points elsewhere. Based on the optimality conditions of Sec. \ref{sec:interp}, consistent change point recovery can be interpreted as the possibility of choosing the initial value of $m_t$ and its slopes so that the graph of $z_t$ lies within $-\lambda$ and $\lambda$, touching the boundaries only within an $O(\epsilon)$ of the true change point instants.

To get some intuition behind the change detection consistency result, notice that a bias term of order $\mu$ in the slope of $m_t$ may lead to a bias of order $\mu (N / M)^2$ between the end points of $z_t$ in one segment, so for a given $\lambda$, $\mu$ has to be of order $\lambda M^2 / N^2$ to achieve the interpolation constraints. Therefore, by choosing $\lambda = o(N^2)$, trend filtering can choose the slopes within $\mu$ of the true slopes of $m_t^o$ so that $z_t$ touches alternating boundaries of the tube $\pm \lambda$ within an $O(\epsilon)$ neighborhood of the true change points. On the other hand, making $\lambda$ grow slowly with $N$ may allow $z_t$ to touch the boundaries $\pm \lambda$ outside the $O(\epsilon)$ neighborhoods of the true change points, due to the variability of the integrated random walk, leading to ``fake'' change points. This can be prevented by noting that the variance of integrated random walk, around the center of each segment, is of order $\sigma^2 (N/M)^2$, \ie, its standard deviation is of order $\sigma N / M$. Therefore, $\lambda$ should grow faster that $N$ to keep the boundaries away from the random variations of the integrated random walk. Notice, finally, that it is not possible to recover the exact location of the change points, but only approximately, because in the neighborhood of the true change points the graph of $z_t$ stays very close to the boundary, and noise may inevitably introduce fake change points in those neighborhoods (as $z_t$ tries to cross the boundary).

This heuristic description can be formalized, as done in \cite{Rojas-Wahlberg-14b} for the fused lasso, to establish that for $\lambda \propto N^c$, with $1 < c < 2$, $\ell_1$ trend filtering achieves approximate ($O(\epsilon)$) change point recovery with probability tending to $1$ as $N \to \infty$, if all consecutive changes in slope have alternating signs.

In case \emph{some of the consecutive changes in slope have the same sign}, change point consistency is not possible. This follows again from the dual interpretation of $\ell_1$ trend filtering. When two consecutive change points have the same direction, $z_t$ is forced to go from $\lambda$ ($-\lambda$) to $\lambda$ ($-\lambda$) within the segment joining the change points, without ever crossing the boundary in between. Due to these interpolation constraints, the deterministic term in $m_t - y_t$ is asymptotically negligible, so $z_t$ must stay very close to the boundary without crossing it within the segment (outside the $O(\epsilon)$ neighborhoods of the true change points); due to the random component of $m_t - y_t$, the probability of achieving this does not go to zero as $N \to \infty$, leading to the possible appearance of fake change points in such segment. This issue is related to the so-called \emph{stair-case effect} in the fused lasso \cite{Rojas-Wahlberg-14b}, where the presence of two or more consecutive changes of the mean level in the same direction introduces spurious change points. An example of this phenomenon will be given in Section~\ref{sec:example}.

\section{A Modified Scheme for Trend Filtering} \label{sec:monitoring}
The discussion in Sec. \ref{sec:interp} leads to a natural scheme for achieving change point consistency even in the presence of consecutive changes in slope of the same sign. The key idea is that, asymptotically in $N$, fake change points can only appear in segments between other detected change points. Therefore, if we apply $\ell_1$ trend filtering to a $N$-sample sequence, the first and last detected change points are \emph{real}, \ie, they approximately correspond to true change points. We can then consider only the segment of data between the first and last change points, and apply $\ell_1$ trend filtering to this new data. Proceeding iteratively in this manner, we can single out all the true change points of the sample, disregarding those fake ones that appear in the first iterations of this scheme.


\section{Examples} \label{sec:example}

\begin{figure}[h]
\centering
\includegraphics[width=0.7\columnwidth]{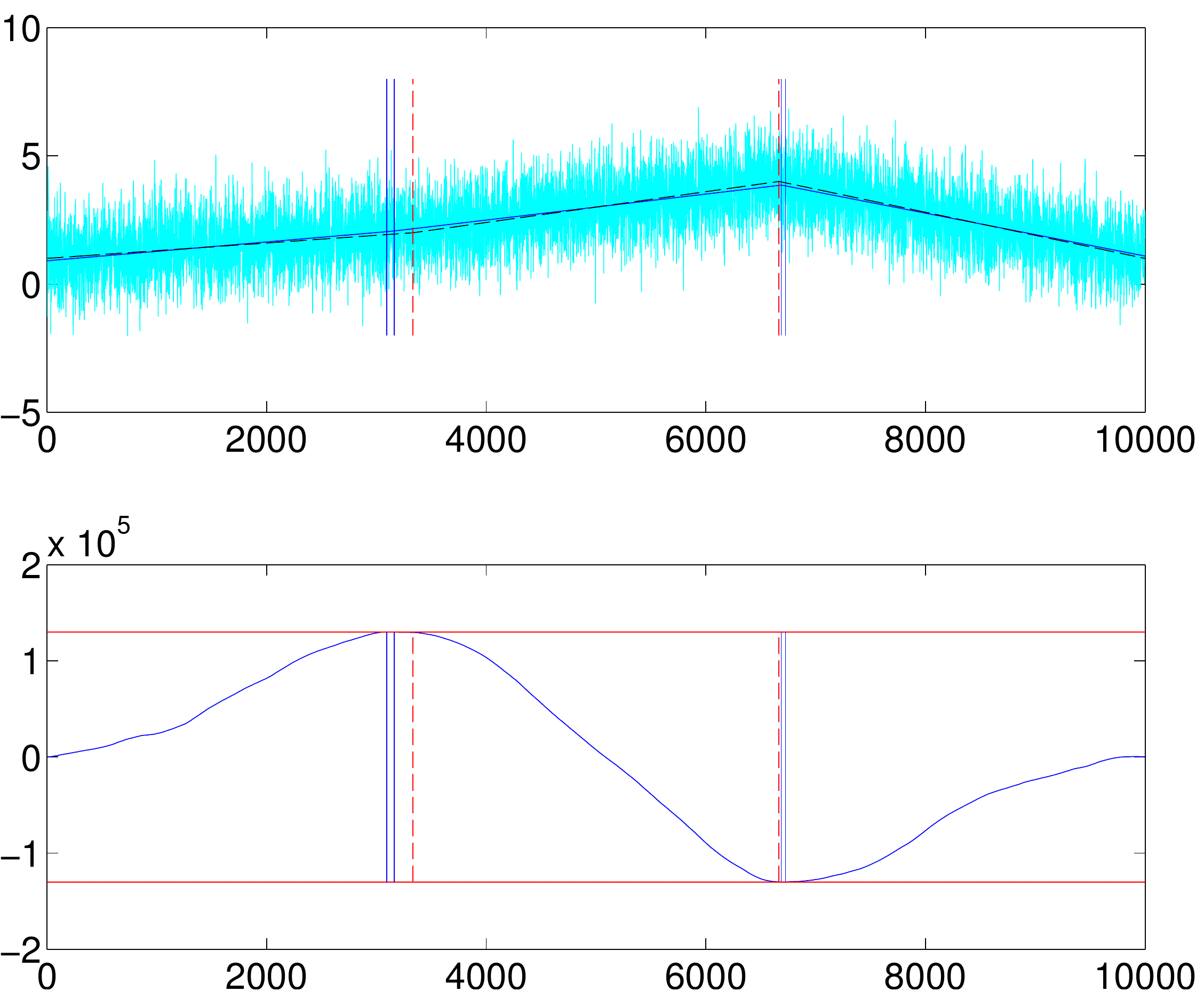}
\caption{Successful change point recovery. \textbf{Top:} The cyan line shows the data $y_t$, while the dashed (black) and solid (blue) lines correspond to the true and estimated means, respectively, which nearly coincide; the dashed (red) and solid (blue) vertical lines denote the true and estimated change points, respectively. \textbf{Bottom:} Graph of the dual variable $z_t$.}
\label{fig:2}
\end{figure}

In this section we consider two examples, both with $N = 10000$ samples and variance $\sigma^2 = 1$. In the first example, the mean value is a piecewise linear signal which goes from $1$ to $2$ between $t=1$ and $t = 3333$, then to $4$ at $t= 6666$, and finally back to $1$ at $t=10000$. The slopes of this signal change in alternating directions, so from Sec.~\ref{sec:interp} we should expect $\ell_1$ trend filtering to achieve change detection recovery; the situation is shown in Fig.~\ref{fig:2}. Here, $\ell_1$ trend filtering successfully detects change points in the neighborhood of their true locations, and no spurious change points have appeared. Actually, 2 estimated change points appear close to the first true one, but due to their close proximity we consider them as one successful detection. $\lambda$ was chosen equal to $130000$, which is approximately $N^{1.3}$; for this dataset, the method can recover detect the change points when $\lambda \leqslant 260000 \approx N^{1.35}$.

\begin{figure}[h]
\centering
\includegraphics[width=0.7\columnwidth]{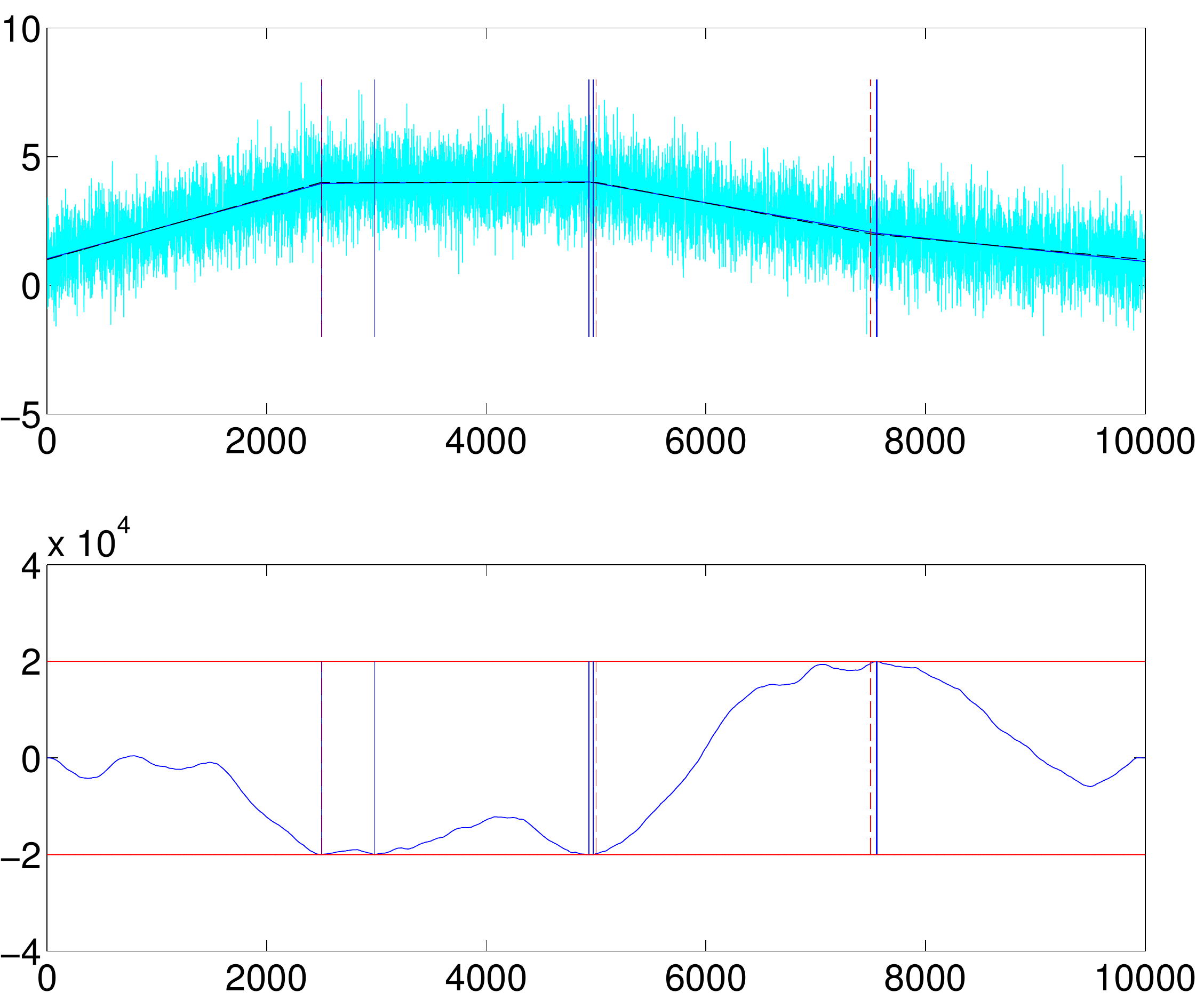}
\caption{Failed change point recovery. Same notation as Fig.~\ref{fig:2}.}
\label{fig:3}
\end{figure}

Consider now a second example, where the mean is a piecewise linear signal which goes from $1$ to $4$ between $t=1$ and $t = 2500$, stays at $4$ until $t= 5000$, then decreases down to $2$ at $t= 7500$, and finally goes back to $1$ at $t=10000$. In this case, the changes in slope are not purely alternating in sign, so we should expect the presence of fake change points not close to the true ones. Fig.~\ref{fig:3} shows this situation for $\lambda = 20000$. Here we see that $\ell_1$ trend filtering correctly detects the true change points (\ie, it identifies change points close to the true ones); however, there is a fictitious change point in the segment between the first two true change points. Notice that changing $\lambda$ has no effect on this fake change point: it is not possible to remove it by increasing $\lambda$, as this cannot alter the bias term in the affected segment, but only on those segments where $z_t$ is forced to move from one boundary to the other (or close to the initial and final end-points).

\begin{figure}[h]
\centering
\includegraphics[width=0.7\columnwidth]{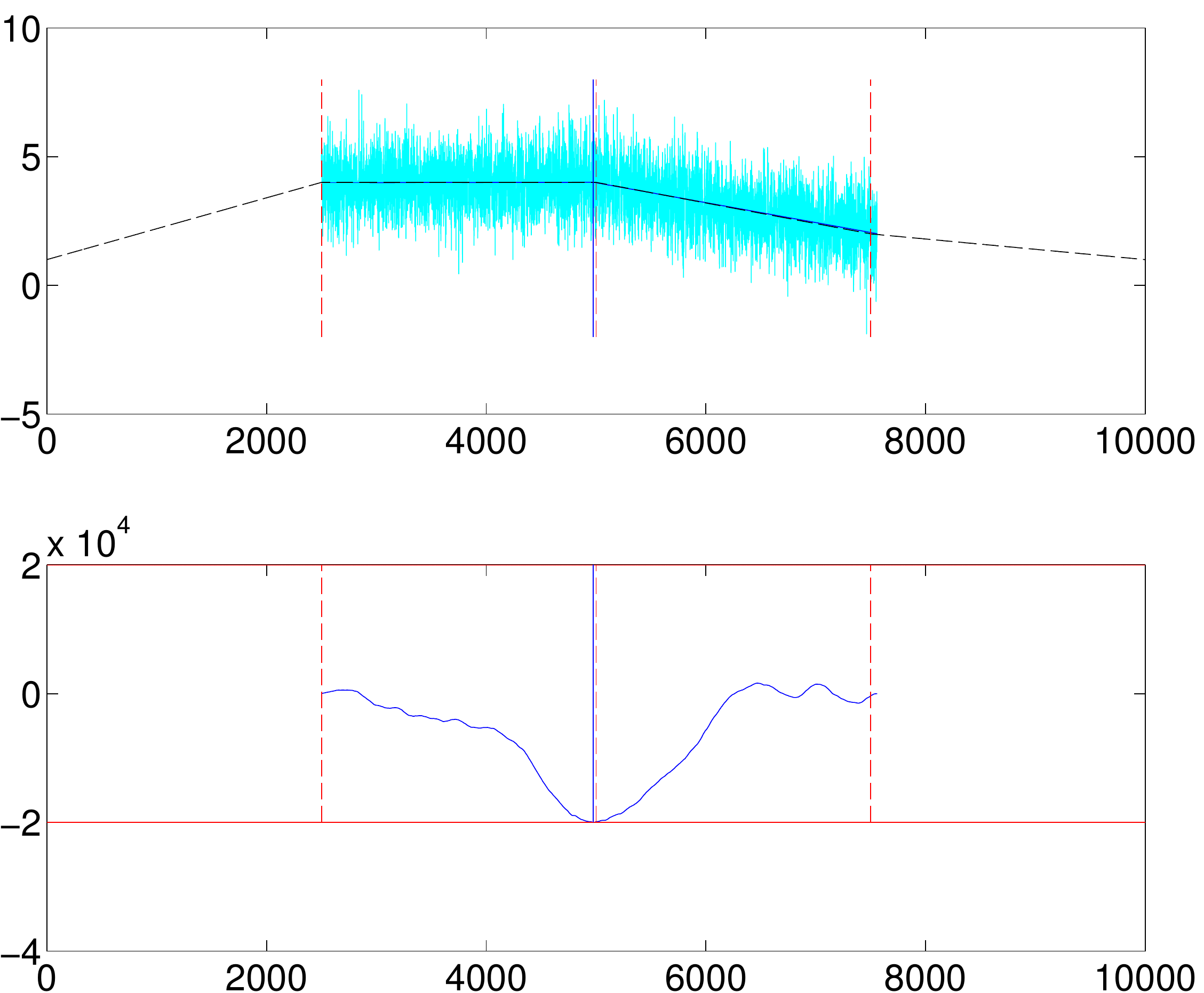}
\caption{Modified scheme applied to second example.}
\label{fig:4}
\end{figure}

To remove the presence of fake change points, we use the scheme of Sec.~\ref{sec:monitoring}, according to which we consider the first and last estimated change points as ``true'' ones, and re-apply $\ell_1$ trend filtering only to the data between them. The result is shown in Fig.~\ref{fig:4}. Note here that the fake change point has completely disappeared! Furthermore, since $z_t$ at the location of the fake change point is far from $\pm \lambda$, the monitoring scheme provides a very robust means to remove such artifice.

\section{Conclusions} \label{sec:conclusion}

In this paper we have studied the change point consistency of $\ell_1$ trend filtering, a technique for the estimation of piecewise linear trends in noisy time series. Based on a geometric interpretation of the method, we have provided an intuitive understanding of situations when the method succeeds and when it fails. Furthermore, building on this interpretation, we have developed a technique for removing false change points.

\vfill
\pagebreak

\bibliographystyle{IEEEbib}
\bibliography{cristian,refs_bo}

\begin{thebibliography}{10}

\bibitem{gustafsson-2000a}
F.~Gustafsson,
\newblock {\em Adaptive Filtering and Change Detection},
\newblock John Wiley \& Sons, 1 edition, Sept. 2000.

\bibitem{Kim-Koh-Boyd-Gorinevsky-09}
S.-J. Kim, K.~Koh, S.~Boyd, and D.~Gorinevsky,
\newblock ``$l_1$ trend filtering,''
\newblock {\em SIAM Review}, vol. 51(2), pp. 339--360, 2009.

\bibitem{citeulike:11033787}
R.~J. Hodrick and E.~C. Prescott,
\newblock ``{Postwar U.S. Business Cycles: An Empirical Investigation},''
\newblock {\em Journal of Money, Credit and Banking}, vol. 29, no. 1, pp. 1+,
  Feb. 1997.

\bibitem{OhlssonLB:10}
H.~Ohlsson, L.~Ljung, and S.~Boyd,
\newblock ``Segmentation of {ARX}-models using sum-of-norms regularization,''
\newblock {\em Automatica}, vol. 46, pp. 1107 -- 1111, 2010.

\bibitem{Ryan11}
R.~Tibshirani and J.~Taylor,
\newblock ``The solution path of the generalized lasso,''
\newblock {\em Annals of Statistics}, vol. 39, no. 3, pp. 1335--1371, 2011.

\bibitem{SDN06}
G.~Steidl, S.~Didas, and J.~Neumann,
\newblock ``Splines in higher order tv regularization,''
\newblock {\em International Journal of Computer Vision}, vol. 70, pp.
  241--255, 2006.

\bibitem{Ryan14}
R.~Tibshirani,
\newblock ``Adaptive piecewise polynomial estimation via trend filtering,''
\newblock {\em Annals of Statistics}, vol. 42, no. 1, pp. 285--323, 2014.

\bibitem{citeulike:9394683}
M.~A. Little and N.~S. Jones,
\newblock ``{Generalized methods and solvers for noise removal from piecewise
  constant signals. I. Background theory},''
\newblock {\em Proceedings of the Royal Society A: Mathematical, Physical and
  Engineering Science}, vol. 467, no. 2135, pp. 3088--3114, 2011.

\bibitem{citeulike:9394684}
M.~A. Little and N.~S. Jones,
\newblock ``{Generalized methods and solvers for noise removal from piecewise
  constant signals. II. New methods},''
\newblock {\em Proceedings of the Royal Society A: Mathematical, Physical and
  Engineering Science}, vol. 467, no. 2135, pp. 3115--3140, 2011.

\bibitem{Rinaldo2009}
A.~Rinaldo,
\newblock ``Properties and refinements of the fused lasso,''
\newblock {\em The Annals of Statistics}, vol. 37, no. 5B, pp. pp. 2922--2952,
  2009.

\bibitem{Rojas-Wahlberg-14b}
C.~R. Rojas and B.~Wahlberg,
\newblock ``On change point detection using the fused lasso method,''
\newblock {\em Annals of Statistics (submitted for publication)}, 2014,
\newblock arXiv:1401.5408.

\bibitem{Piotr14}
P.~Fryzlewicz,
\newblock ``Wild binary segmentation for multiple change-point detection,''
\newblock {\em Annals of Statistics},
\newblock to appear.

\end{thebibliography}

\end{document}